\documentclass[12pt]{article}

\usepackage{latexsym, amsfonts, amsmath, amssymb}
\newcommand{\be}{\begin{equation}}
\newcommand{\ee}{\end{equation}}
\newcommand{\bea}{\begin{eqnarray}}
\newcommand{\eea}{\end{eqnarray}}
\newcommand{\bean}{\begin{eqnarray*}}
\newcommand{\eean}{\end{eqnarray*}}

\newcommand{\brray}{\begin{array}}
\newcommand{\erray}{\end{array}}
\newcommand{\ben}{\begin{equation}{nonumber}}
\newcommand{\een}{\end{equation}{nonumber}}

\newtheorem{dfn}{Definition}[section]
\newtheorem{thm}[dfn]{Theorem}
\newtheorem{lmma}[dfn]{Lemma}
\newtheorem{ppsn}[dfn]{Proposition}
\newtheorem{crlre}[dfn]{Corollary}
\newtheorem{xmpl}[dfn]{Example}
\newtheorem{rmrk}[dfn]{Remark}
\newcommand{\bdfn}{\begin{dfn}}
\newcommand{\bthm}{\begin{thm}}
\newcommand{\blmma}{\begin{lmma}}
\newcommand{\bppsn}{\begin{ppsn}}
\newcommand{\bcrlre}{\begin{crlre}}
\newcommand{\bxmpl}{\begin{xmpl}}
\newcommand{\brmrk}{\begin{rmrk}}
\newcommand{\edfn}{\end{dfn}}
\newcommand{\ethm}{\end{thm}}
\newcommand{\elmma}{\end{lmma}}
\newcommand{\eppsn}{\end{ppsn}}
\newcommand{\ecrlre}{\end{crlre}}
\newcommand{\exmpl}{\end{xmpl}}
\newcommand{\ermrk}{\end{rmrk}}

\newcommand{\clh}{{\cal H}}

\newcommand{\clq}{{\cal Q}}

\def\a*{{\cal A}_{h,*}}
\def\B{{\cal B}(h)}
\def\B1{{\cal B}_1(h)}
\def\b{{\cal B}^{\rm s.a.}(h)}
\def\b1{{\cal B}^{\rm s.a.}_1(h)}

\newcommand{\ot}{\otimes}

\newcommand{\raro}{\rightarrow}

\def \qed {$\Box$}

\def\a*{{\cal A}_{h,*}}
\def\B{{\cal B}(h)}
\def\B1{{\cal B}_1(h)}
\def\b{{\cal B}^{\rm s.a.}(h)}
\def\b1{{\cal B}^{\rm s.a.}_1(h)}

\begin{document}

\begin{center}
{\Large{\bf Example of a group whose quantum isometry group does not depend on the generating set }}
\end{center}
\begin{center}
{\large { Arnab Mandal}}
\end{center} 
\begin{abstract} 
In this article we have shown that the quantum isometry group of $C_r^*(\mathbb{Z})$, denoted by $\mathbb{Q}(\mathbb{Z},S)$ as in \cite{gos_man}, with respect to a symmetric generating set $S$ does not depend on the generating set $S$. Moreover, we have proved that the result is no longer true if the group $\mathbb{Z}$ is replaced by   $\underbrace{\mathbb{Z} \times \mathbb{Z} \times\cdots \times \mathbb{Z}}_{n \ copies}$ for $n>1$.
\end{abstract}

\section{Introduction}
Quantum groups are very important mathematical entities which appear in several areas of mathematics and physics, often as same kind of generalised symmetry objects. Beginning from the pioneering work by Drinfeld, Jimbo, Manin, Woronowicz and others nearly three decades ago (\cite{Drinfeld} \cite{Manin}, \cite{wor} and references therein) there is now a vast literature on quantum groups both from algebraic and analytic (operator algebraic) viewpoints. Generalizing group actions on spaces, notions of (co)actions of quantum groups on possibly noncommutative spaces have been formulated and studied by many mathematicians. In this context, S.Wang \cite{wang} came up with the definition of quantum automorphism groups of certain mathematical structures (typically finite sets, matrix algebras etc.) and such quantum groups have been studied in depth since then. Later on, a number of mathematicians including Banica, Bichon and others (\cite{ban_1} \cite{finite graph} and references therein) developed a theory of quantum automorphism groups of finite metric spaces and finite graphs.  With a more geometric setup in 2009, Goswami \cite{Gos} defined and proved existence of an analogue of the group of isometries of a Riemannian manifold, in the framework of the so-called compact quantum groups \`a la Woronowicz. 
In fact, he considered the more general setting of noncommutative manifold, given by spectral triples defined by Connes \cite{con} and under some mild regularity conditions, he proved the existence of a universal compact quantum group (termed as the quantum isometry group) acting on the $C^*$-algebra underlying the noncommutative manifold such that the action also commutes with a natural analogue of Laplacian of the spectral triple. Furthermore, Goswami and Bhowmick formulated \cite{qorient}  the notion of a quantum group analogue of the group of orientation preserving isometries and its existence as the universal object in a suitable category was proved.  After that, several authors studied quantum isometry groups of different spectral triples in recent years.\\
In literature, we have an interesting as well as important spectral triple on $C_r^*(\Gamma)$ \cite{connes_compact}, coming from the word length of a finitely generated discrete group $\Gamma$ corresponding to a symmetric generating set, say S. There have been several articles already on computations and study of the quantum isometry groups of such spectral triples, e.g. (\cite{org filt}, \cite{free cyclic}, \cite{grp algebra}, \cite{gos_man}, \cite{mandal}, \cite{S_n}, \cite{doubling}, \cite{dihedral} and references therein). We denote it by $\mathbb{Q}(\Gamma,S)$ as in \cite{gos_man}. It is already known that in general $\mathbb{Q}(\Gamma,S_1)$ and $\mathbb{Q}(\Gamma,S_2)$ are not isomorphic for different choices of $S_1$ and $S_2$. Indeed, they are drastically different for certain choices of generating sets. For instance, if we choose n such that $g.c.d(n,4)=1$, then the group $\mathbb{Z}_n \times \mathbb{Z}_4$ is isomorphic to $\mathbb{Z}_{4n}$. Consider the generating sets $S_1=\lbrace(1,0),(0,1), (-1,0), (0,-1) \rbrace$ and $S_2=\lbrace (1,1), (-1,-1) \rbrace$ respectively for $\mathbb{Z}_n \times \mathbb{Z}_4 $. The underlying $C^*$-algebra of $\mathbb{Q}(\mathbb{Z}_n \times \mathbb{Z}_4,S_1)$ is noncommutative by Theorem 4.10 of \cite{gos_man}. On the other hand, $\mathbb{Q}(\mathbb{Z}_n \times \mathbb{Z}_4,S_2)$ is the doubling of $C^*(\mathbb{Z}_n \times \mathbb{Z}_4)$ corresponding to the automorphism given by $a \mapsto a^{-1} \ \forall \ a \in \mathbb{Z}_n \times \mathbb{Z}_4$ from \cite{grp algebra}. Hence its underlying $C^*$-algebra is commutative. They are non-isomorphic even in the vector space level. In this context, it is quite natural to find out the groups whose quantum isometry group does not depend on the generating set. Our main goal of this article is to provide one such example.   \\
The paper is organized as follows. In Section 2 we recall some definitions and necessary facts regarding to compact quantum groups, quantum isometry groups and the doubling procedure of a compact quantum group. Section 3 contains the main results of this paper. In Theorem 3.1 we have proved that the quantum isometry group of $C_r^*(\mathbb{Z})$ remains unchanged if we change the generating sets. Theorem \ref{new thm} tells us that this  is no longer true if $\mathbb{Z}$ is replaced by   $\underbrace{\mathbb{Z} \times \mathbb{Z} \times\cdots \times \mathbb{Z}}_{n \ copies}$ for $n>1$. 
\section{Preliminaries}
First of all, we fix some notational convention. The algebraic tensor product and spatial (minimal) $C^*$-tensor product are denoted  by $\ot$ and $\hat{\ot}$  respectively throughout the article. We'll use the leg-numbering notation. Let $\mathcal{H}$ be a complex Hilbert space, $\mathcal{K}(\mathcal{H})$ the $C^*$-algebra of compact operators on it, and $\mathcal{Q}$ a unital $C^*$-algebra. The multiplier algebra $\mathcal{M}(\mathcal{K}(\mathcal{H})\hat{\ot} \mathcal{Q})$ has two natural embeddings into $\mathcal{M}(\mathcal{K}(\mathcal{H})\hat{\ot} \mathcal{Q} \hat{\ot} \mathcal{Q})$, one obtained by extending the map $x \mapsto x \ot 1$ and the second one is obtained by composing this map with the flip on the last two factors. We will write $\omega^{12}$ and $\omega^{13}$ for the images of an element $\omega \in \mathcal{
M}(\mathcal{K}(\mathcal{H})\hat{\ot} \mathcal{Q})$ under these two maps respectively. We'll denote by $\mathcal{H} \bar{\ot} \mathcal{Q}$ the Hilbert $C^*$-module obtained by completing $\mathcal{H} \ot \mathcal{Q}$ with respect to the norm induced by the $\mathcal{Q}$ valued inner product $\langle\langle \xi \ot q, \xi^{\prime} \ot q^{\prime}\rangle\rangle \ := \langle\xi,\xi^{\prime}\rangle q^*q^{\prime}$, where $\xi,\xi^{\prime} \in \mathcal{H}$ and $q,q^{\prime} \in \mathcal{Q}$.

\subsection{Compact quantum groups}
 In this subsection we recall some standard definitions related to compact quantum groups. We recommend \cite{wor}, \cite{Van} for more details.
\bdfn
A compact quantum group (CQG in short) is a pair $(\clq,\Delta)$, where $\clq$ is a unital $C^{\ast}$- algebra and $\Delta:\clq\raro\clq\hat{\ot}\clq$ is a unital $\ast$-homomorphism (called the comultiplication), such that
\begin{enumerate}
\item $(\Delta\ot id)\Delta=(id\ot \Delta)\Delta$ as homomorphism $\clq\raro\clq\hat{\ot}\clq\hat{\ot}\clq$ (coassociativity).
\item  The spaces $\Delta(\clq)(1\ot \clq)= Span \{\Delta(b)(1\ot a)|a,b\in\clq\}$ and $(1\ot \clq)\Delta(\clq)$ are dense in $\clq\hat{\ot}\clq$.
\end{enumerate}
\edfn
Sometimes we may denote the CQG $(\mathcal{Q},\Delta)$ simply as $\mathcal{Q}$, if $\Delta$ is clear from the context.
\bdfn
A CQG morphism from $(\mathcal{Q}_1,\Delta_1)$ to another $(\mathcal{Q}_2,\Delta_2)$ is a unital $C^*$-homomorphism $\pi: \mathcal{Q}_1 \mapsto \mathcal{Q}_2$ such that $(\pi \ot \pi)\Delta_1=\Delta_2 \pi$.
\edfn

\bdfn
 We say that a CQG $(\mathcal{Q},\Delta)$ acts on a unital $C^*$-algebra $\mathcal{B}$ if there is a unital $C^*$-homomorphism (called action) $\alpha : \mathcal{B} \rightarrow \mathcal{B} \hat{\ot} \mathcal{Q} $ satisfying the following:
 \begin{enumerate}
 \item $(\alpha \ot id)\alpha = (id \ot \Delta)\alpha $.
 \item Linear span of $\alpha(\mathcal{B})(1 \ot \mathcal{Q})$ is norm-dense in $\mathcal{B} \hat{\ot} \mathcal{Q}$.
 \end{enumerate}
\edfn

\bdfn
Let $(\clq,\Delta)$ be a CQG. A unitary representation of $\clq$ on a Hilbert
space $\clh$ is a $\mathbb{C}$-linear map $U$ from $\clh$ to the Hilbert module
$\clh\bar{\ot}
\clq$ such that
\begin{enumerate}
\item $\langle\langle U(\xi),U(\eta)\rangle\rangle=\langle \xi,\eta \rangle 1_{\clq}$, where $\xi,\eta\in \clh$.
\item  $(U\ot {\rm id})U=(id\ot \Delta)U$.
\item $Sp \ \{U(\clh)\clq\}$ is dense in $\clh\bar{\ot}\clq$.
\end{enumerate}
\edfn
Given such a unitary representation we have a unitary element $\tilde{U}$ belonging to $\mathcal{M}(\mathcal{K}(\mathcal{H})\hat{\ot} \mathcal{Q})$ given by $\tilde{U}(\xi \ot b)= U(\xi)b,(\xi \in \mathcal{H}, \ b \in \mathcal{Q})$ satisfying $(id \ot \Delta)(\tilde{U})=\tilde{U}^{12}\tilde{U}^{13}$.

\subsection{Quantum isometry groups}
In \cite{Gos} Goswami introduced the notion of quantum isometry group of a spectral triple satisfying certain regularity conditions. We refer to \cite{Gos}, \cite{qorient}, \cite{org filt} for the original formulation of quantum isometry groups and its various avatars including the quantum isometry group for orthogonal filtrations.
\bdfn
Let $(\mathcal{A}^{\infty},\mathcal{H},\mathcal{D})$ be a spectral triple of compact type (a la Connes). Consider the category $Q(\mathcal{D})\equiv Q(\mathcal{A}^{\infty},\mathcal{H},\mathcal{D})$ whose objects are $(\mathcal{Q}, U)$, where $(\mathcal{Q},\Delta)$ is a CQG having a unitary representation U on the Hilbert space $\mathcal{H}$ satisfying the following:
\begin{enumerate}
\item $\tilde{U}$ commutes with $(\mathcal{D}\ot 1_{\mathcal{Q}})$.
\item $(id \ot \phi)\circ ad_{\tilde{U}}(a) \in (\mathcal{A}^{\infty})^{\prime\prime}$ for all $a \in \mathcal{A}^{\infty}$ and $\phi$ is any state on $\mathcal{Q}$, where $ad_{\tilde{U}}(x): = \tilde{U}(x \ot 1)\tilde{U}^*$ for $x \in \mathcal{B}(\mathcal{H})$. 
 \end{enumerate}
A morphism between two such objects  $(\mathcal{Q}, U)$ and $ (\mathcal{Q}^\prime, U^\prime)$ is a CQG morphism $\psi : \mathcal{Q} \rightarrow \mathcal{Q}^\prime $ such that $U^\prime = (id \ot \psi)U$. If a universal object exists in $Q(\mathcal{D})$ then we denote it by $\widetilde{QISO^{+}(\mathcal{A}^{\infty},\mathcal{H},\mathcal{D})}$ and the corresponding largest Woronowicz subalgebra for which $ad_{\tilde{U_0}}$ is faithful, where $U_{0}$ is the unitary representation of $\widetilde{QISO^{+}(\mathcal{A}^{\infty},\mathcal{H},\mathcal{D})}$, is called the quantum group of orientation preserving isometries and denoted by $QISO^{+}(\mathcal{A}^{\infty},\mathcal{H},\mathcal{D})$.
\edfn
Let us state Theorem $2.23$ of \cite{qorient} which gives a sufficient condition for the existence of  $QISO^{+}(\mathcal{A}^{\infty},\mathcal{H},\mathcal{D})$.
\bthm\label{existence thm}
 Let $(\mathcal{A}^{\infty},\mathcal{H},\mathcal{D})$ be a spectral triple of compact type. Assume that $\mathcal{D}$ has one dimensional kernel spanned by a vector $\xi \in \mathcal{H} $ which is cyclic and separating for $\mathcal{A}^{\infty}$ and each eigenvector of \ $\mathcal{D}$ belongs to $\mathcal{A}^{\infty}\xi$. Then QISO$^+(\mathcal{A}^{\infty},\mathcal{H},\mathcal{D})$ exists. 
 \ethm
 Here we briefly discuss a specific case of interest for us.  For more details see Section 2.2 of \cite{gos_man}. Let $\Gamma$ be a finitely generated discrete group with a symmetric  generating set  $S$ not containing the identitity of $\Gamma$
 (symmetric means $g \in S$ if and only if $g^{-1} \in S$)
 and let $l$ be the corresponding word length function. We  define an operator $D_{\Gamma}$ by $D_\Gamma(\delta_g)=l(g)\delta_g$,
 where $\delta_g$ denotes the vector in $l^{2}(\Gamma)$ which takes value $1$ at the point $g$ and $0$ at all other points. Observe that $\delta_g, g \in \Gamma$ forms 
  an orthonormal basis of $l^2(\Gamma)$. Let $\tau$ be the canonical trace on the group $C^*$-algebra given by $\tau(\sum_{g\in \Gamma} c_g\lambda_g)=c_e$, where e is the identity element of the group $\Gamma$. Connes first considered this spectral triple $(\mathbb{C}\Gamma$, $l^2(\Gamma),D_\Gamma)$ in \cite{connes_compact}. It is easy to check that $ (\mathbb{C}\Gamma$, $l^2(\Gamma),D_\Gamma)$ is a spectral triple using Lemma 1.1 of \cite{oza_rie}. Moreover,  QISO$^+(\mathbb{C}\Gamma$, $l^2(\Gamma),D_\Gamma)$ exists by Theorem \ref{existence thm}, 
   taking $\delta_e$ as the cyclic separating vector for $ \mathbb{C}\Gamma$. Note that its action $\alpha$ (say) on $C_r^*(\Gamma)$ is determined by 
   $$ \alpha(\lambda_{\gamma})= \sum_{\gamma^{\prime} \in S} \lambda_{\gamma^{\prime}} \ot q_{\gamma, \gamma^{\prime}},$$
   where the matrix $[q_{\gamma,\gamma^{\prime}}]_{\gamma, \gamma^{\prime} \in S} \in M_{card(S)}(\mathbb{Q}(\Gamma,S))$ is called the fundamental unitary of $\mathbb{Q}(\Gamma,S)$.

\subsection{Doubling of a CQG}\label{doub section}
We briefly recall the doubling procedure of a compact quantum group from \cite{S_n}, \cite{doubling}.
Let $(\mathcal{Q},\Delta)$ be a CQG with a CQG-automorphism $\theta$ such that $\theta^2=id$. The doubling of this CQG, 
 say $(\mathcal{D}_{\theta}(\mathcal{Q}),\tilde{\Delta})$, is  given by  $\mathcal{D}_{\theta}(\mathcal{Q}) :=\mathcal{Q}\oplus \mathcal{Q}$ (direct sum as a $C^*$-algebra),
 and the coproduct is defined by the following, where we have denoted  the injections of $\mathcal{Q}$ onto 
 the first and second coordinate in $\mathcal{D}_{\theta}(\mathcal{Q})$ by $\xi$ and $\eta$ respectively, i.e. 
$\xi(a)=(a,0), \  \eta(a)= (0,a), \ (a \in \mathcal{Q}).$
$$\tilde{\Delta} \circ \xi= (\xi \ot \xi + \eta \ot [\eta \circ \theta])\circ  \Delta,$$
$$\tilde{\Delta} \circ \eta= (\xi \ot \eta + \eta \ot [\xi \circ \theta])\circ \Delta.$$
 
 \section{Main Results}
 Before going to the main theorem we make one convention. Inverse of any element $x \in \mathbb{Z}$ is denoted by $-x$. We will also follow the same convention for  $\underbrace{\mathbb{Z} \times \mathbb{Z} \times\cdots \times \mathbb{Z}}_{n \ copies},$ where $ n>1$. 
 \bthm\label{main thm}
 For any symmetric generating set $S$, the quantum isometry group $\mathbb{Q}(\mathbb{Z},S)$ is isomorphic to $\mathcal{D}_{\theta}(C^*(\mathbb{Z}))$ with respect to the automorphism $\theta$ given by $\theta(x)=-x \ \forall \ x \in \mathbb{Z}$.  
 \ethm
  {\it Proof :}\\
 Let us assume that $S$ is any generating set for $\mathbb{Z}$, i.e. $S=\lbrace a_1, -a_1, \cdots, a_k, -a_k \rbrace $. Without loss of generality we can assume that $a_i > 0 \ \forall \ i$ and $a_1 < a_2 < \cdots < a_k$. We would like to mention here that the largest number $a_k$ and the smallest number $-a_k$ of the generating set $S$ will play a crucial role in the proof. For each $i=1,\cdots,k-1$ there exists positive integers $c_i, d_i $ such that $c_i a_i= d_i a_k \ \forall \ i=1,\cdots, k-1$. Moreover, $c_i > d_i$ as $a_i < a_k \forall \ i=1,\cdots,k-1$.  
 Now the action $\alpha$ of $\mathbb{Q}(\mathbb{Z},S)$ on $C_r^*(\mathbb{Z})$ is defined as \\
\begin{eqnarray*}
 \alpha(\lambda_{a_{1}}) &=& \lambda_{a_{1}} \ot A_{11} + \lambda_{-a_{1}} \ot A_{12} + \lambda_{a_{2}} \ot A_{13} + \lambda_{-a_{2}} \ot A_{14} + \cdot \cdot \cdot     + \\ 
 &&\lambda_{a_{k}} \ot A_{1(2k-1)} + \lambda _{-a_{k}} \ot A_{1(2k)}, \\
  \alpha(\lambda_{-a_{1}}) &=& \lambda_{a_{1}} \ot A_{12}^* + \lambda_{-a_{1}} \ot A_{11}^* + \lambda_{a_{2}} \ot A_{14}^* + \lambda_{-a_{2}} \ot A_{13}^* + \cdot \cdot \cdot     + \\ 
  &&\lambda_{a_{k}} \ot A_{1(2k)}^* + \lambda _{-a_{k}} \ot A_{1(2k-1)}^*, \\
 \alpha(\lambda_{a_{2}}) &=& \lambda_{a_{1}} \ot A_{21} + \lambda_{-a_{1}} \ot A_{22} + \lambda_{a_{2}} \ot A_{23} + \lambda_{-a_{2}} \ot A_{24} + \cdot \cdot \cdot     + \\
&& \lambda_{a_{k}} \ot A_{2(2k-1)} + \lambda _{-a_{k}} \ot A_{2(2k)}, \\
\alpha(\lambda_{-a_{2}}) &=& \lambda_{a_{1}} \ot A_{22}^* + \lambda_{-a_{1}} \ot A_{21}^* + \lambda_{a_{2}} \ot A_{24}^* + \lambda_{-a_{2}} \ot A_{23}^* + \cdot \cdot \cdot     + \\ 
&&\lambda_{a_{k}} \ot A_{2(2k)}^* + \lambda _{-a_{k}} \ot A_{2(2k-1)}^*, \\
\vdots  &&   \hspace{3cm}           \vdots  \hspace{4cm}           \vdots\\
\vdots  &&   \hspace{3cm}           \vdots  \hspace{4cm}           \vdots\\
\vdots  &&   \hspace{3cm}           \vdots  \hspace{4cm}           \vdots\\
\vdots  &&   \hspace{3cm}           \vdots  \hspace{4cm}           \vdots\\
\alpha(\lambda_{a_{k}}) &=& \lambda_{a_{1}} \ot A_{k1} + \lambda_{-a_{1}} \ot A_{k2} + \lambda_{a_{2}} \ot A_{k3} + \lambda_{-a_{2}} \ot A_{k4} + \cdot \cdot \cdot     + \\ 
 &&\lambda_{a_{k}} \ot A_{k(2k-1)} + \lambda _{-a_{k}} \ot A_{k(2k)}, \\
\alpha(\lambda_{-a_{k}}) &=& \lambda_{a_{1}} \ot A_{k2}^* + \lambda_{-a_{1}} \ot A_{k1}^* + \lambda_{a_{2}} \ot A_{k4}^* + \lambda_{-a_{2}} \ot A_{k3}^* + \cdot \cdot \cdot     + \\ 
  &&\lambda_{a_{k}} \ot A_{k(2k)}^* + \lambda _{-a_{k}} \ot A_{k(2k-1)}^*. \\ 
\end{eqnarray*} 
 The fundamental unitary is of the form 
$$ U=
 \begin{pmatrix}
 A_{11} & A_{12} & A_{13} & A_{14} & \cdots & A_{1(2k-1)} & A_{1(2k)}\\
 A_{12}^* & A_{11}^* & A_{14}^* & A_{13}^* & \cdots &  A_{1(2k)}^* & A_{1(2k-1)}^*\\
 A_{21} & A_{22} & A_{23} & A_{24} & \cdots & A_{2(2k-1)} & A_{2(2k)}\\
 A_{22}^* & A_{21}^* & A_{24}^* & A_{23}^* & \cdots & A_{2(2k)}^* & A_{2(2k-1)}^*\\
 \vdots  && \hspace{1cm}        \vdots \\
 A_{k1} & A_{k2} & A_{k3} & A_{k4} & \cdots & A_{k(2k-1)} & A_{k(2k)}\\
 A_{k2}^* & A_{k1}^* & A_{k4}^* & A_{k3}^* & \cdots &  A_{k(2k)}^* & A_{k(2k-1)}^*\\ 
 \end{pmatrix}.$$\\
 
 Our aim is to show that it reduces to the following: \\\\
\begin{equation}\label{eqn 1} 
 \begin{pmatrix}
 A_{11} & A_{12} & 0 & 0 & \cdots & 0 & 0\\
 A_{12}^* & A_{11}^* & 0 & 0 & \cdots &  0 & 0\\
 0 & 0 & A_{23} & A_{24} & \cdots & 0 & 0\\
 0 & 0 & A_{24}^* & A_{23}^* & \cdots & 0 & 0\\
 \vdots  && \hspace{1cm}        \vdots \\
 0 & 0 & 0 & 0 & \cdots & A_{k(2k-1)} & A_{k(2k)}\\
 0 & 0 & 0 & 0 & \cdots &  A_{k(2k)}^* & A_{k(2k-1)}^*\\ 
 \end{pmatrix}
 \end{equation}
i.e. only the diagonal $(2\times2)$ block survives and others become zero. First we will show that it reduces to the form \\
 $$\begin{pmatrix}
 A_{11} & A_{12} & A_{13} & A_{14} & \cdots & 0 & 0\\
 A_{12}^* & A_{11}^* & A_{14}^* & A_{13}^* & \cdots &  0 & 0\\
 A_{21} & A_{22} & A_{23} & A_{24} & \cdots & 0 & 0\\
 A_{22}^* & A_{21}^* & A_{24}^* & A_{23}^* & \cdots & 0 & 0\\
 \vdots  && \hspace{1cm}        \vdots \\
 0 & 0 & 0 & 0 & \cdots & A_{k(2k-1)} & A_{k(2k)}\\
 0 & 0 & 0 & 0 & \cdots &  A_{k(2k)}^* & A_{k(2k-1)}^*\\ 
 \end{pmatrix}$$\\ 
 
i.e. $A_{i(2k-1)}=0, A_{i(2k)}=0 \ \forall \ i=1,\cdots,k-1$. Using the antipode $A_{k(2j-1)}=A_{k(2j)}=0 \ \forall \ j=1,\cdots, k-1$. We break the proof into a number of lemmas. 
\blmma\label{lem1}
$A_{i(2k-1)}^{c_i}= A_{i(2k)}^{c_i}=0 \ \forall \ i=1,\cdots, k-1.$  
\elmma 
 {\it Proof :}\\
 Consider the term $\alpha(\lambda_{c_ia_i})= \alpha(\lambda_{d_ia_k})$ forall $i=1,\cdots,k-1$. Comparing the coefficients of $\lambda_{c_ia_k}$ and $\lambda_{c_i(-a_k)}$ on both sides of the relation $\alpha(\lambda_{c_ia_i})= \alpha(\lambda_{d_ia_k})$ we obtain $A_{i(2k-1)}^{c_i}=A_{i(2k)}^{c_i}=0$ as the right hand side of the equation does not contain any terms with coefficients $\lambda_{c_ia_k}$ and $\lambda_{c_i(-a_k)}$ as well. \qed\\
 
 Our goal is to show that $A_{i(2k-1)}$ and $A_{i(2k)}$ are normal $\forall \ i=1,\cdots,k-1$. Then by Lemma \ref{lem1} one can conclude $A_{i(2k-1)}=A_{i(2k)}=0$.
 
\blmma \label{lem3}
 If $a_p + a_q=a_l + a_m$ for some $p,q,l,m \in \mathbb{N}$, then $A_{k(2p-1)}A_{k(2q-1)}= A_{k(2l-1)}A_{k(2m-1)}.$
 \elmma 
 {\it Proof :}\\
Using the relation $\alpha(\lambda_{a_p + a_q})=\alpha(\lambda_{a_l + a_m})$ comparing the coefficient of $\lambda_{2a_k}$ on both sides we have $A_{p(2k-1)}A_{q(2k-1)}=A_{l(2k-1)}A_{m(2k-1)}$. Applying the antipode we obtain  $A_{k(2q-1)}^*A_{k(2p-1)}^*=A_{k(2m-1)}^*A_{k(2l-1)}^*$. This implies $A_{k(2p-1)}A_{k(2q-1)}$ $= A_{k(2l-1)}A_{k(2m-1)}.$ \qed\\\\
We state three auxiliary lemmas (Lemma \ref{lem4} to Lemma \ref{lem6}) whose proof will follow by exactly the same arguments used in Lemma \ref{lem3}. We omit the proofs.   

\blmma\label{lem4}
If $a_p + a_q = -a_l + a_m$, then $A_{k(2p-1)}A_{k(2q-1)}= A_{k(2l)}A_{k(2m-1)}.$ 
\elmma

\blmma\label{lem5}
If $a_p - a_q = -a_l+a_m$, then $A_{k(2p-1)}A_{k(2q)}= A_{k(2l)}A_{k(2m-1)}.$ 
\elmma

\blmma\label{lem6}
If $a_p - a_q=-a_l-a_m$, then $A_{k(2p-1)}A_{k(2q)}= A_{k(2l)}A_{k(2m)}.$ 
\elmma 

\blmma\label{lem7}
$A_{k(2i-1)}A_{k(2i)}=A_{k(2i)}A_{k(2i-1)}=0 \ \forall \ i=1,\cdots, k$.
\elmma
{\it Proof :}\\
Comparing the coefficients of $\lambda_{2a_k}$ and $\lambda_{-2a_k}$ from the relation $\alpha(\lambda_{a_i}). \alpha(\lambda_{-a_i})= \lambda_e \ot 1$ one can get $A_{i(2k-1)}A_{i(2k)}^*=A_{i(2k)}A_{i(2k-1)}^*=0$. Applying the antipode we have $A_{k(2i)}^*A_{k(2i-1)}^*=0$ which implies $A_{k(2i-1)}A_{k(2i)}=0$. Similarly we can get $A_{k(2i)}A_{k(2i-1)}=0$ from the relation $\alpha(\lambda_{-a_i}). \alpha(\lambda_{a_i})= \lambda_e \ot 1$. \qed\\

\blmma\label{lem8}
$A_{ki}A_{kj}=0 \ \forall \ i,j$ with $i \neq j$.
\elmma     
{\it Proof :}\\
 We will show that $A_{k(2j)}A_{k1}=A_{k(2j-1)}A_{k1}=0$ forall $j=2,\cdots,k$. By Lemma \ref{lem7} we have $A_{k2}A_{k1}=0$. Then $A_{ki}A_{k1}=0$ forall $i$ with $i \neq 1$ will be proved. Other relations will follow by repeating the same line of arguments. Suppose for some fixed $j \neq 1$, we have $a_1+ a_j=a_l-a_m=a_g+a_t=2a_s$ for some $l,m,g,t,s \in \mathbb{N}$, where the collection  $\lbrace g,t \rbrace$ is considered  with $g \neq t$ and  $g,t \in \lbrace 2,3,\cdots, j-1, j+1, \cdots,k \rbrace$. Also assume if $\lbrace g,t \rbrace$ is considered then we will not consider $\lbrace t,g \rbrace$.  Sometimes the sets $\lbrace l,m \rbrace, \lbrace g,t \rbrace$ and $\lbrace s \rbrace$  may be empty depending on the choice of the generating set $S$. Note that the collections $\lbrace l,m \rbrace$ and $\lbrace g,t \rbrace$ are not necessarily singleton but finite. Now from the condition $\alpha(\lambda_{a_k}). \alpha(\lambda_{-a_k})= \lambda_e \ot 1$ comparing the coefficient of $\lambda_{a_1+a_j}$ one can deduce
 \begin{eqnarray}\label{eqn 2}
 A_{k1}A_{k(2j)}^* + A_{k(2j-1)}A_{k2}^* + \sum_{g,t} [ A_{k(2g-1)}A_{k(2t)}^*  + A_{k(2t-1)}A_{k(2g)}^*] \nonumber \\ + \sum_{l,m} [ A_{k(2l-1)}A_{k(2m-1)}^* + A_{k(2m)}A_{k(2l)}^*] +A_{k(2s-1)}A_{k(2s)}^*= 0 
\end{eqnarray}
Multiplying $A_{k(2j)}$ and $A_{k1}^*$ on the left hand side and right hand side respectively of the equation (\ref{eqn 2}) we get 
\begin{eqnarray}\label{eqn 3}
 A_{k(2j)}A_{k1}A_{k(2j)}^*A_{k1}^* + A_{k(2j)}A_{k(2j-1)}A_{k2}^*A_{k1}^* + \sum_{g,t} A_{k(2j)}A_{k(2g-1)}A_{k(2t)}^*A_{k1}^* \nonumber \\ + \sum_{g,t} A_{k(2j)}A_{k(2t-1)}A_{k(2g)}^*A_{k1}^* + \sum_{l,m} A_{k(2j)}A_{k(2l-1)}A_{k(2m-1)}^*A_{k1}^* \nonumber \\ + \sum_{l,m} A_{k(2j)}A_{k(2m)}A_{k(2l)}^*A_{k1}^*+ A_{k(2j)}A_{k(2s-1)}A_{k(2s)}^*A_{k1}^*=0
\end{eqnarray}
Now the relation $a_1+a_j=a_j+a_1$ gives us $a_1-a_j=-a_j+a_1$. By Lemma \ref{lem5} we have $A_{k1}A_{k(2j)}=A_{k(2j)}A_{k1}$. Similarly, $A_{k(2j)}A_{k(2g-1)}= A_{k1}A_{k(2t)}$ as $a_1-a_t= -a_j+a_g$ from the assumed condition $a_1+a_j= a_g+a_t$. Moreover, by Lemma \ref{lem4}, Lemma \ref{lem5} and Lemma \ref{lem6} we get 
 $$A_{k(2j)}A_{k(2t-1)}= A_{k1}A_{k(2g)}, A_{k(2j)}A_{k(2l-1)}=A_{k1}A_{k(2m-1)},$$  
 $$A_{k(2j)}A_{k(2s-1)}= A_{k1}A_{k(2s)}, A_{k(2j)}A_{k(2m)}=A_{k1}A_{k(2l)},$$
 as $-a_j+a_t= a_1-a_g, -a_j+a_l=a_1+a_m, -a_j +a_s=a_1-a_s$ and $ -a_j-a_m= a_1-a_l$ respectively. Using these relations and Lemma \ref{lem7} the equation (\ref{eqn 3}) reduces to 
\begin{eqnarray}\label{eqn 4}
A_{k(2j)}A_{k1}(A_{k(2j)}A_{k1})^* + \sum_g A_{k(2j)}A_{k(2g-1)}(A_{k(2j)}A_{k(2g-1)})^* + \nonumber \\  \sum_t A_{k(2j)}A_{k(2t-1)}(A_{k(2j)}A_{k(2t-1)})^* + \sum_l A_{k(2j)}A_{k(2l-1)}(A_{k(2j)}A_{k(2l-1)})^* \nonumber \\ + \sum_m A_{k(2j)}A_{k(2m)}(A_{k(2j)}A_{k(2m)})^* + A_{k(2j)}A_{k(2s-1)}( A_{k(2j)}A_{k(2s-1)})^*=0 
\end{eqnarray} 
This shows that $A_{k(2j)}A_{k1}=0$ as the left hand side of the equation (\ref{eqn 4}) is the sum of positive elements of a $C^*$-algebra. Similarly from the relation $\alpha(\lambda_{a_k}). \alpha(\lambda_{-a_k})= \lambda_e \ot 1$ comparing the coefficient of $\lambda_{a_1-a_j}$ we get that $A_{k(2j-1)}A_{k1}=0$.     
     \qed\\  
\blmma \label{lem9}
$A_{i(2k)}$ and $A_{i(2k-1)}$ are normal $\forall \ i=1,\cdots,k$.
\elmma  
 {\it Proof :}\\
 By the unitarity condition of $U$ we have 
 $$ \sum_i A_{k(2i-1)}A_{k(2i-1)}^* + \sum_i A_{k(2i)}A_{k(2i)}^*=1,  $$
 $$ \sum_i A_{k(2i-1)}^*A_{k(2i-1)} + \sum_i A_{k(2i)}^*A_{k(2i)}=1.  $$
 Thus $A_{k(2i-1)}^2 A_{k(2i-1)}^*= A_{k(2i-1)}$ and $A_{k(2i-1)}^* A_{k(2i-1)}^2= A_{k(2i-1)}$ by using Lemma \ref{lem8}. Hence, $A_{k(2i-1)}A_{k(2i-1)}^*= A_{k(2i-1)}^*A_{k(2i-1)}^2 A_{k(2i-1)}^*=A_{k(2i-1)}^*A_{k(2i-1)}$. Applying the antipode we get that $A_{i(2k-1)}$ is normal. Similarly, it can be shown that $A_{i(2k)}$ is normal. \qed\\\\
 By Lemma \ref{lem1} and Lemma \ref{lem9} we have $A_{i(2k)}=A_{i(2k-1)}=0 \ \forall \ i=1,\cdots, k-1.$ Repeating the same arguments using from Lemma \ref{lem1} to Lemma \ref{lem9} we can conclude that the fundamental unitary finally reduces to the form as in (\ref{eqn 1}), i.e.\\
 $$ \begin{pmatrix}
 A_{11} & A_{12} & 0 & 0 & \cdots & 0 & 0\\
 A_{12}^* & A_{11}^* & 0 & 0 & \cdots &  0 & 0\\
 0 & 0 & A_{23} & A_{24} & \cdots & 0 & 0\\
 0 & 0 & A_{24}^* & A_{23}^* & \cdots & 0 & 0\\
 \vdots  && \hspace{1cm}        \vdots \\
 0 & 0 & 0 & 0 & \cdots & A_{k(2k-1)} & A_{k(2k)}\\
 0 & 0 & 0 & 0 & \cdots &  A_{k(2k)}^* & A_{k(2k-1)}^*\\ 
 \end{pmatrix}.$$\\
 Note that $A_{i(2i-1)}A_{i(2i)}=A_{i(2i)}A_{i(2i-1)}=0 \ \forall \ i=1,\cdots,k$ and all the entries of the fundamental unitary are normal. Moreover for each $i$,  $A_{i(2i-1)}A_{i(2i-1)}^*$ and $A_{i(2i)}A_{i(2i)}^*$ are projections. We also have $A_{i(2i-1)}^2 A_{i(2i-1)}^*= A_{i(2i-1)}$ and $A_{i(2i)}^2 A_{i(2i)}^*= A_{i(2i)}$.  For every $a_i$ and $a_j$, there exists positive integers $p(i,j)$ and $q(i,j)$ depending on $i,j$ such that $p(i,j)a_i=q(i,j)a_j$. Using this condition one can easily get $A_{i(2i)}^{p(i,j)}=A_{j(2j)}^{q(i,j)}$ and $A_{i(2i-1)}^{p(i,j)}=A_{j(2j-1)}^{q(i,j)} \ \forall \ i,j=1,\cdots,k$. This gives us $A_{i(2i-1)}A_{j(2j)}^{q(i,j)}(A_{j(2j)}^{q(i,j)})^*= 0$ as $A_{i(2i-1)}A_{i(2i)}=0$. Thus $A_{i(2i-1)}A_{j(2j)}A_{j(2j)}^*= 0$ by using that $A_{j(2j)}A_{j(2j)}^*$ is projection and $A_{j(2j)}$ is normal. Finally we get that $A_{i(2i-1)}A_{j(2j)}=0$ as $A_{j(2j)}^2 A_{j(2j)}^*= A_{j(2j)}$ and $A_{j(2j)}$ is normal. Similarly one can deduce that $A_{i(2i)}A_{j(2j-1)}=0 \ \forall \ i,j=1,\cdots,k$.\\\\ 
We can define the map from $\mathbb{Q}(\mathbb{Z},S)$ to $\mathcal{D}_{\theta}(C^*(\mathbb{Z}))$ by
$$ A_{i(2i-1)} \mapsto (\lambda_{a_i},0),$$
$$ A_{i(2i)} \mapsto (0, \lambda_{-a_i})=(0, \lambda_{\theta(a_i)}),$$
 $ \forall \ i=1,\cdots,k.$ Clearly this gives an isomorphism between two CQG's. \qed.\\
\brmrk
From now on, the quantum isometry group $\mathbb{Q}(\mathbb{Z},S)$ can be written simply as $\mathbb{Q}(\mathbb{Z})$. 
\ermrk

Now we are going to show that the quantum isometry group of $\underbrace{\mathbb{Z} \times \mathbb{Z} \times\cdots \times \mathbb{Z}}_{n \ copies}$ for $n>1$ depends on the generating set. We present the case $n=2$ for the simplicity of the exposition. The proof for any $n$ can be adapted similarly from the proof of Theorem \ref{new thm}. Let $S^{\prime}=\lbrace  (1,0),(0,1), (-1,0),(0,-1)\rbrace$ and $S^{\prime\prime}=\lbrace  (1,0),(0,1), (-1,0),(0,-1), (2,0), (-2,0)\rbrace$  be the two different generating sets for $\mathbb{Z} \times \mathbb{Z}$.   
 \bthm \label{new thm}
 $\mathbb{Q}(\mathbb{Z} \times \mathbb{Z}, S^{\prime})$ and $\mathbb{Q}(\mathbb{Z} \times \mathbb{Z}, S^{\prime\prime})$ are not isomorphic to each other.
 \ethm
{\it Proof :}\\
First of all, note that by Proposition 2.29 and Theorem 4.1 of \cite{gos_man} we get that  $\mathbb{Q}(\mathbb{Z} \times \mathbb{Z}, S^{\prime}) \cong C((\mathbb{T} \times \mathbb{T}) \rtimes (\mathbb{Z}_2 ^{2} \rtimes S_2))$. We will show that  
$\mathbb{Q}(\mathbb{Z} \times \mathbb{Z}, S^{\prime\prime})$ is different from  $C((\mathbb{T} \times \mathbb{T}) \rtimes (\mathbb{Z}_2 ^{2} \rtimes S_2))$. Let us assume that $a=(1,0), b=(0,1)$ and $c=(2,0)$. The action $\alpha$ is defined by \\
$$\alpha(\lambda_{a})= \lambda_{a} \ot A + \lambda_{-a} \ot B + \lambda_{b} \ot C + \lambda_{-b} \ot D + \lambda_{c} \ot E + \lambda_{-c} \ot F, $$
$$\alpha(\lambda_{-a})= \lambda_{a} \ot B^* + \lambda_{-a} \ot A^* + \lambda_{b} \ot D^* + \lambda_{-b} \ot C^* + \lambda_{c} \ot F^* + \lambda_{-c} \ot E^*, $$
$$\alpha(\lambda_{b})= \lambda_{a} \ot G + \lambda_{-a} \ot H + \lambda_{b} \ot I + \lambda_{-b} \ot J + \lambda_{c} \ot K + \lambda_{-c} \ot L, $$
$$\alpha(\lambda_{-b})= \lambda_{a} \ot H^* + \lambda_{-a} \ot G^* + \lambda_{b} \ot J^* + \lambda_{-b} \ot I^* + \lambda_{c} \ot L^* + \lambda_{-c} \ot K^*, $$
$$\alpha(\lambda_{c})= \lambda_{a} \ot M + \lambda_{-a} \ot N + \lambda_{b} \ot O + \lambda_{-b} \ot P + \lambda_{c} \ot Q + \lambda_{-c} \ot R, $$
$$\alpha(\lambda_{-c})= \lambda_{a} \ot N^* + \lambda_{-a} \ot M^* + \lambda_{b} \ot P^* + \lambda_{-b} \ot O^* + \lambda_{c} \ot R^* + \lambda_{-c} \ot Q^*. $$
 The fundamental unitary is of the form 
$$\begin{pmatrix}
A & B & C & D & E & F\\
B^* & A^* & D^* & C^* & F^* & E^*\\
G & H & I & J & K & L\\
H^* & G^* & J^* & I^* & L^* & K^*\\
M & N & O & P & Q & R\\
N^* & M^* & P^* & O^* & R^* & Q^*\\
\end{pmatrix}.$$\\

Note that the product of any two different elements of each row of the fundamental unitary is zero by the similar arguments of Lemma \ref{lem8}. Moreover, using  Lemma \ref{lem9} we get that all the entries are normal. Using the relation $\alpha(\lambda_{2a})=\alpha(\lambda_c)$ comparing the coefficients of $\lambda_b$ and $\lambda_{-b}$ on both sides we have $O=P=0$. Applying the antipode and involution we get $K=L=0$. Similarly, comparing the coefficients of $\lambda_{2b}, \lambda_{-2b},\lambda_{2c}$ and $\lambda_{-2c}$ from the same condition we have $C^2=D^2=E^2=F^2=0$ as well. This gives us $C=D=E=F=0$ as they are normal. Using the antipode and the involution we obtain $G=H=M=N=0$. \\ 

Thus the fundamental unitary reduces to
$$\begin{pmatrix}
A & B & 0 & 0 & 0 & 0\\
B^* & A^* & 0 & 0 & 0 & 0\\
0 & 0 & I & J & 0 & 0\\
0 & 0 & J^* & I^* & 0 & 0\\
0 & 0 & 0 & 0 & Q & R\\
0 & 0 & 0 & 0 & R^* & Q^*\\
\end{pmatrix}.$$\\ 
Observe that $A^2=Q, B^2=R$ comparing the coefficients of $\lambda_c$ and $\lambda_{-c}$ from the condition $\alpha(\lambda_{2a})=\alpha(\lambda_c)$. Moreover, $AR=BQ=0$ as $AB=QR=0$. The underlying $C^*$-algebra of $\mathbb{Q}(\mathbb{Z} \times \mathbb{Z}, S^{\prime\prime})$ is generated by the elements $A,B,I$ and $J$. We also get that $AI=IA, AJ=JA, BI=IB$ and $BJ=JB$ comparing the coefficients of $\lambda_{a+b}, \lambda_{-a+b}, \lambda_{a-b}$ and $\lambda_{-a-b}$    on both sides from the relation $\alpha(\lambda_{a+b})=\alpha(\lambda_{b+a})$. Clearly, the CQG $\mathbb{Q}(\mathbb{Z} \times \mathbb{Z}, S^{\prime\prime})$ is identified with $\mathbb{Q}(\mathbb{Z}) \hat{\ot} \mathbb{Q}(\mathbb{Z})$.  Note that the underlying $C^*$-algebra of $\mathbb{Q}(\mathbb{Z}) \hat{\ot} \mathbb{Q}(\mathbb{Z})$ is isomorphic to $[C^*(\mathbb{Z}) \oplus C^*(\mathbb{Z})] \hat{\ot} [C^*(\mathbb{Z}) \oplus C^*(\mathbb{Z})]$. The isomorphism is defined as follows:
$$ A \mapsto  (\lambda_1,0) \ot 1,$$
$$ B \mapsto  (0, \lambda_{-1}) \ot 1,$$
$$ I \mapsto 1 \ot (\lambda_1,0),$$
$$ J \mapsto 1 \ot (0, \lambda_{-1}),$$
where $\lbrace 1,-1 \rbrace$ is the standard minimal generating set for $\mathbb{Z}$. Thus $\mathbb{Q}(\mathbb{Z} \times \mathbb{Z}, S^{\prime\prime}) \cong C((\mathbb{T} \times \mathbb{T}) \rtimes \mathbb{Z}_2 ^{2} ) $ as $\mathbb{Q}(\mathbb{Z})$ is isomorphic to $C(\mathbb{T} \rtimes \mathbb{Z}_2)$. It is clearly not isomorphic with $C((\mathbb{T} \times \mathbb{T}) \rtimes (\mathbb{Z}_2 ^{2} \rtimes S_2))$, hence we are done. \qed

 E-mail address: {\bf arnabmaths@gmail.com}\\\\
 SCHOOL OF MATHEMATICAL SCIENCES, NATIONAL INSTITUTE OF SCIENCE EDUCATION AND RESEARCH BHUBANESWAR, HBNI, JATNI, 752050, INDIA.

\end{document}